\documentclass[a4paper,11pt]{article}

\usepackage{latexsym}
\usepackage{amsmath,amsthm,amsfonts,amssymb}
\usepackage{graphicx,color,epsfig}
\usepackage{nicefrac}

\def\R{\mathbb{R}}
\def\N{\mathbb{N}}
\def\Z{\mathbb{Z}}

\def\leq{\leqslant}
\def\geq{\geqslant}

\newtheorem{Pro}{Proposition}[section]
\newtheorem{Lem}[Pro]{Lemma}
\newtheorem{Thm}[Pro]{Theorem}

\begin{document}

\title{Extensive escape rate in lattices of weakly coupled expanding maps with holes}
\author{Jean-Baptiste Bardet$^\ast$ and Bastien Fernandez$^\dag$}
\date{}   
\maketitle
\begin{center}
$^\ast$ Laboratoire de Math\'ematiques Rapha\" el Salem\footnote{UMR 6085 CNRS - Universit\'e de Rouen} \\ Avenue de l'Universit\'e\\ 76801 Saint \'Etienne du Rouvray, France

$^\dag$ Centre de Physique Th\'{e}orique\footnote{UMR 6207 CNRS - Universit\'e Aix-Marseille II - Universit\'e Aix-Marseille I - Universit\'e Sud Toulon-Var}\\ CNRS Luminy Case 907\\ 13288 Marseille CEDEX 9, France
\end{center}

\begin{abstract}
This paper discusses possible approaches to the escape rate in infinite lattices of weakly coupled maps with uniformly expanding repeller. It is proved that computed-via-volume rates of spatially periodic approximations grow linearly with the period size, suggesting normalized escape rate as the appropriate notion for the infinite system. The proof relies on symbolic dynamics and is based on the control of cumulative effects of perturbations within cylinder sets. A piecewise affine diffusive example is presented that exhibits monotonic decay of the escape rate with coupling intensity. 

\end{abstract}

\maketitle

\section{Introduction}
Coupled Map Lattices (CML) have been introduced in the early 1980's as space-time discrete models of reaction-diffusion processes and of other spatially extended systems \cite{CF05,K93}. Originally designed to ascertain the stability of numerical simulation schemes, CML soon appeared to provide a unique opportunity to extend the theory of dynamical systems to realistic infinite dimensional examples with discrete time. In particular, a large effort has been made to define and prove the existence and uniqueness of the (infinite-dimensional analogue of) SRB measure for small perturbations of uniformly expanding maps, see \cite{KL06} and the expanded list of references therein. 

Of special interest in dynamical systems with spatial extension is also the definition of dynamical quantifiers such as the entropy \cite{ACFM02,AMU04,CE99} or the escape rate. 
The escape rate comes as a natural estimator in presence of 'holes', {\sl i.e.}\ regions in phase space where the dynamics is not defined or where the orbit is lost. Such types of systems often occur in concrete examples; billiards, lattice or hard-spheres gases, etc \cite{GD95}. When the generator of the dynamics is expanding, the existence of holes suggests to consider the repeller, {\sl i.e.}\ the invariant set of points for which the orbit never reaches any hole. A natural problem is then to characterize the exponential rate of escape from a small neighborhood of the repeller. 

The goal of this paper is to suggest a suitable definition of the escape rate for infinite lattices of weakly coupled maps with holes. The theory of (finite dimensional) dynamical system contains two distinct approaches to the escape rate, either via the conditionally invariant measure \cite{PY79} or by volume estimates \cite{BR75,ER85}. None of these approaches has an immediate extension to the infinite dimensional setting. On one hand, no proof of existence of an absolutely continuous conditionally invariant measure for CML is available in the literature, even in the simplest setting. On the other hand, the second approach needs to be adapted in order to ensure finiteness of estimated (Lebesgue) volumes at all times. 

For simplicity, we deal here with CML with convolution couplings and transitive expanding real maps satisfying the Pianigiani-Yorke conditions \cite{PY79} and we consider the perturbative regime where the symbolic dynamics is well under control (section 2). We first discuss the existence of conditionally invariant measure and some properties of the related escape rate for spatially periodic approximations (section 3.1). A piecewise affine example is analyzed in section 3.2 which shows linear growth of the escape rate with the size of the approximation and also monotonic decay as the coupling intensity increases. This behavior motivates an alternative approach to escape rate that deals with volume estimates of points remaining in the neighborhood of the repeller (section 3.3). In this context, we prove the existence of the escape rate for every spatially periodic approximation of arbitrary systems and its linear growth with the system size. The corresponding asymptotic slope can be regarded as an escape rate per lattice unit for the infinite system. The proof is given in section 4 and essentially relies on the control of cumulative effects of perturbations within cylinder sets.

\section{Weakly coupled lattices of Pianigiani-Yorke maps}
A CML on the lattice $\Z$ is a dynamical system generated by a self-map $F$ in $\ell^\infty(\Z)$ (endowed with uniform topology and associated norm $\|\cdot\|$) whose specificity resides in the following expression. The mapping $F$ is given by the composition $F=C\circ F_0$ of the direct product $F_0$ of identical real maps and of a coupling $C$ \cite{CF05,K93}. In simplest situations, the coupling is the {\bf convolution operator} associated with an arbitrary non-negative and normalized sequence $\{c_n\}_{n\in\Z}$ in $\ell^1(\Z)$, {\sl i.e.}\ we have 
\[
(Cx)_s=\sum_{n\in\Z} c_n x_{s-n},\ \forall s\in\Z,
\]
for every configuration $\{x_s\}_{s\in\Z}\in \ell^\infty(\Z)$ (see \cite{AF00} for properties of convolution operators in $\ell^\infty(\Z)$). The individual mapping $F_0$ explicitly writes
\[
(F_0 x)_s=f(x_s),\ \forall s\in\Z.
\]
where $f:\R\mapsto\R$ is a (locally bounded) map  that verifies the following conditions. There exists a finite collection $\{\text{I}_i\}_{i\in {\mathcal A}}$ (where ${\mathcal A}=\{1,\cdots,N\}$, $N\in\N$, $N>1$) of pairwise disjoint bounded and closed intervals with the properties that for each $i\in {\mathcal A}$  
\begin{itemize}
\item[(i)] $f$ is of class $C^2$ on $\text{I}_i$ and $\inf\limits_{\text{I}_i}\left|f'\right|>1$,
\item[(ii)] there exists $j\in \mathcal{A}$ such that $\text{I}_j\subset \text{Int}\; f(\text{I}_i)$. Moreover, for every $j$ such that $f(\text{I}_i)\cap \text{I}_j\neq \emptyset$, we have $\text{Int}\; f(\text{I}_i)\supset \text{I}_j$.
\item[(iii)] In addition, we also assume that $f$ is {\bf transitive}, {\sl i.e.}\ for every pair $i,j\in \mathcal{A}$ there is a $t\in\N$ such that $\text{I}_j\subset f^t(\text{I}_i)$.
\end{itemize}
We will denote $\text{I}=\bigcup_i\text{I}_i$. In this context, the region $f(I)\setminus I$ plays the role of holes for the dynamical system $(f,I)$. 

\noindent
The conditions on $f$ imply several properties of the dynamics, see e.g.\ \cite{CMS94,CMT98,PY79}, in particular
\begin{itemize}
\item[$\bullet$] the structural stability with respect to $C^1$ perturbations, 
\item[$\bullet$] the existence of a repeller $\kappa_f$ on which the dynamics is conjugated to a topological Markov chain ($\kappa_f$ is a Cantor set, forward invariant under $f$ and defined by $\kappa_f=\bigcap\limits_{t\in\N}f^{-t}(\text{I})$), 
\item[$\bullet$] the existence of a conditionally invariant probability measure that is absolutely continuous with respect to the Lebesgue measure.
\end{itemize}
The topological Markov chain operates in the set $\Omega\subset {\mathcal A}^\N$  of admissible sequences which are constrained by the condition (ii) ({\sl i.e.}\ the sequence $\{\vartheta^t\}_{t\in\N}$ is {\bf admissible} iff every pair of consecutive symbols $\vartheta^t,\vartheta^{t+1}$ complies with the condition $\text{I}_{\vartheta^{t+1}}\subset \text{Int}\; f(\text{I}_{\vartheta^t})$). Each point in $\kappa_f$ is defined by the intersection
\[
\bigcap_{t\in\N}f^{-t}\left(\text{I}_{\vartheta^t}\right)
\]
for some associated code $\{\vartheta^t\}_{t\in\N}\in \Omega$. 

All these properties immediately translate to the uncoupled map lattice $F_0$ (obtained for $C=\text{Id}$, the identity mapping and provided that one only considers periodic lattices when the existence of quasi-invariant measure is concerned). In particular, the symbolic dynamics is given by the action of time translations on sequences of symbol configurations in $\Omega^\Z$. Introducing the notation $\text{\bf I}_\theta=\bigotimes\limits_{s\in\Z}\text{I}_{\theta_s}$ where $\theta=\{\theta_s\}_{s\in\Z}$ is an arbitrary {\bf symbol configuration} in ${\mathcal A}^\Z$, the repeller $K_0=\kappa_f^\Z$ can be regarded as the set of points defined by 
\[
\bigcap_{t\in\N}F_0^{-t}\text{\bf I}_{\theta^t}
\]
for some sequence $\{\theta^t\}_{t\in\N}\in \Omega^\Z$. Alternatively, we have $K_0=\bigcap\limits_{t\in\N}F_{0}^{-t}\text{\bf I}$ where $\text{\bf I}=\text{I}^{\Z}$.

\noindent
Intuitively, the dynamical properties listed above are expected to extend to weakly coupled map lattices with symbolic dynamics being unaffected. To formally state this result, we use the norm on $\mathcal{L}(\ell^\infty(\Z))$, also denoted by $\|\cdot\|$, induced by the uniform norm in $\ell^\infty(\Z)$.  
\begin{Thm} {\rm \cite{AF00}\footnote{In \cite{AF00}, there are two points in the proof of this theorem that are incorrect. One is the proof of a technical statement (namely Proposition 3.2) and the other is the argument that guarantees that cylinder sets are not empty. Both points are however valid and we provide a corrected proof in Appendix \ref{A-ERRATUM}.}}
Given $f$ with properties {\rm (i)--(iii)}, there exists $\epsilon_f>0$ such that for every $C$ so that $\|\text{Id}-C\|<\epsilon_f$ and every sequence $\{\theta^t\}_{t\in\N}\in\Omega^\Z$, the intersection
\[
\bigcap_{t\in\N}F^{-t}\text{\bf I}_{\theta^t}
\]
is non empty and reduces to a point. All such points form a $F$-invariant, closed, and totally disconnected\footnote{but not compact for the uniform topology} set $K=\bigcap\limits_{t\in\N}F^{-t}\text{\bf I}$ on which $F$ is topologically conjugated to the time translation acting on $\Omega^\Z$.
\label{STABSTRUCT}
\end{Thm}
Notice that, due to the conditions $c_n\geq 0$ and ${\displaystyle\sum_{n\in\Z}}c_n=1$, we have $\|\text{Id}-C\|=2(1-c_0)$ meaning that the condition $\|\text{Id}-C\|<\epsilon_f$ actually involves a single entry of the sequence $\{c_n\}$.  

The CML has an obvious {\bf translational symmetry} invariance; namely the map $F$ commutes with translations along the lattice $\Z$. Besides dealing with infinite lattices, this symmetry legitimates focusing on periodic lattices $\Z_L=\Z/L\Z$ (indeed, we have $F\R^{\Z_L}\subset \R^{\Z_L}$ for every $L\in\N$). Restricting to sequences of $L$-periodic symbol configurations in $\Omega^{\Z_L}$  in Theorem \ref{STABSTRUCT} implies that the CML on the periodic lattice of (arbitrary) length $L$ has a limit Cantor set on which the dynamics is conjugated to a topological Markov chain (on a finite alphabet). Following \cite{PY79}, this suggests to investigate the existence of conditionally invariant measures in periodic CML.

\section{Escape rate per lattice unit}
\subsection{Approach via the conditionally invariant measure}
Let $G:X\subset\R^n\to\R^n$ be a mapping acting on some bounded set $X$. A probability measure $\mu$ on $X$ is said to be {\bf conditionally invariant} (or quasi-invariant) if there exists a constant $\alpha\in (0,1)$ such that $\mu\circ G^{-1}=\alpha\mu$ on $X$ \cite{CMS94}. When this measure exists the quantity 
\[
\tilde\gamma=-\log\alpha.
\]
is called the {\bf escape rate}. As mentioned before, the conditions (i)--(iii) guarantee the existence for $f$ of a unique absolutely continuous conditionally invariant measure \cite{PY79}. 

When $\|\text{Id}-C\|<\epsilon_f$, analogous conditions hold for the restriction of the CML to $\R^{\Z_L}$ ($L\in\N$ arbitrary) with the collection $\{\text{\bf I}_\theta\}_{\theta\in {\mathcal A}^{\Z_L}}$ (see Appendix \ref{A-ERRATUM}). It results that, in the weak coupling regime, a unique absolutely continuous conditionally invariant measure exists for the CML on periodic lattices. An explicit example is provided below.

Furthermore, the same conditions (i)--(iii) for $F\big{|}_{\R^{\Z_L}}$ imply the existence of an invariant measure $\nu$ supported on the limit Cantor set that is a Gibbs state for the potential $\log|F\big{|}_{\R^{\Z_L}}'|$ \cite{CMS94}. Accordingly, the following expression for the escape rate applies \cite{CMT98}
\[
\tilde\gamma=\lambda_\nu-h_\nu
\]
where $\lambda_\nu$ is the sum of (positive) Lyapunov exponents and $h_\nu$ is the Kolmog\-or\-ov-Sinai entropy associated with $\nu$. Thus, in principle, one could have access to the escape rate of periodic lattices of weakly coupled maps based on the knowledge of their Lyapunov exponents and metric entropy (or {\sl vice-versa} as in the example below). However, these quantities are not easily accessible in practice, and their behaviour with the spatial period is mostly unknown. Moreover, the problem of existence of an absolutely continuous conditionally invariant measure for the infinite CML, together with a Gibbs state on the corresponding repeller, remains unsolved. 
\begin{figure}[ht]
\centerline{\includegraphics*[height=50mm]{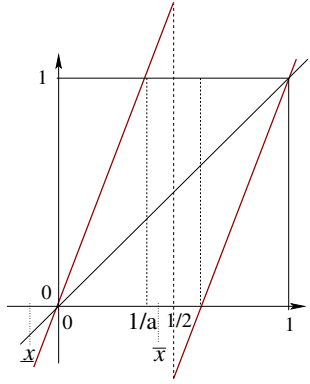}}
\caption{Graph of the Lorenz-type map $f(a)=ax+(1-a)H(x-\frac{1}{2})$ for $a>2$, together with the points $\underline{x}<0<1/a<\overline{x}$.}
\label{LORENZMAP}
\end{figure}

\subsection{Motivating example: CML with Lorenz-type map} 
Consider the coupling operator $C$ derived from the one-dimensional discrete Laplacian
\[
(Cx)_s=x_s+\frac{\epsilon}{2}(x_{s-1}-2x_s+x_{s+1}),\quad \forall s\in\Z
\]
where $\epsilon\in [0,1]$. The local map $f$ is piecewise affine and the points 0 and 1 are fixed points, namely
\[
f(x)=ax+(1-a)H(x-\frac{1}{2}),\quad \forall x\in\R
\]
where $a>2$ and $H$ is the Heaviside function, see Figure \ref{LORENZMAP}. 
This map $f$ clearly satisfies the assumptions (i)--(iii) with $N=2$ and the corresponding intervals are given by 
\[
\text{I}_1=[\underline{x},\overline{x}]\quad\text{and}\quad
\text{I}_2=1-[\underline{x},\overline{x}]
\]
where $\frac{1}{a}<\overline{x}<\frac{1}{2}$ and $\underline{x}<0$ are arbitrary. We have $\text{Int}\; f(\text{I}_i)\supset \text{I}_1\cup \text{I}_2$ for $i=1,2$. 

\noindent
Simple algebra shows that the condition $\|\text{Id}-C\|<\epsilon_f$ in Theorem \ref{STABSTRUCT} becomes $a(1-2\epsilon)>2$ in this case. Indeed, when $a(1-2\epsilon)>2$, for every $\theta\in\{1,2\}^\Z$, the image $F\text{\bf I}_\theta$ covers the product $(\text{I}_1\cup\text{I}_2)^\Z$ (and the expansiveness condition $\|(C\circ F_0)'\|\geq a(1-\epsilon)>1$ holds). This condition is optimal because when $a(1-2\epsilon)\leq 2$, some sequences of symbol configurations in $\{1,2\}^{\Z_2\times\N}$ are not admissible \cite{FG04}. 

Let $L>1$ be fixed in $\N$ and assume that $a(1-2\epsilon)>2$. Then the CML on the $L$-periodic lattice is a piecewise expanding map of constant derivative for which every point in $\left(\text{I}_1\cup\text{I}_2\right)^{\Z_L}$ has $2^L$ pre-images.  
As a consequence, the uniform measure on this set is conditionally invariant. 

\noindent
Computing the determinant of the derivative $(C\circ F_0)'$ yields the following expression for the corresponding escape rate
\begin{equation}
\tilde\gamma_L(\epsilon)=L\log\left(\frac{a}{2}\right)+\sum_{k=1}^L\log\hat{c}\left(\frac{k}{L},\epsilon\right)
\label{TAUXFUITE}
\end{equation}
where $\hat{c}(\omega,\epsilon)=1-\epsilon(1-\cos2\pi\omega)>0$. In particular, since the Lyapunov exponents are in this case trivially given by $\log a\; \hat{c}\left(\frac{k}{L},\epsilon\right)$ ($k=1,\cdots,\lfloor\frac{L}{2}\rfloor$) this expression shows that the Kolmogorov-Sinai entropy of the measure supported on the Cantor set is given by $L\log 2$, {\sl i.e.}\ it corresponds to the topological entropy. 

\noindent
The expression (\ref{TAUXFUITE}) suggests the following additional comments. 
\begin{itemize}
\item[$\bullet$] As a sum of monotonically decreasing functions of $\epsilon$, the escape rate $\tilde\gamma_L(\epsilon)$ decreases when $\epsilon$ increases. In other words, increasing the interaction between sites slows down escape from the repeller, as expected. 
\item[$\bullet$] The rate $\tilde\gamma_L(\epsilon)$ is an extensive quantity that linearly diverges with the period of configurations.  
\end{itemize}
As for the definition of the entropy in spatially extended systems \cite{ACFM02,CE99}, the linear divergence suggests to normalize the escape rate by the period length, {\sl i.e.}\ to consider the ratio
\[
\frac{\tilde\gamma_L(\epsilon)}{L}=\log\left(\frac{a}{2}\right)+\frac{1}{L}\sum\limits_{k=1}^L\log\hat{c}\left(\frac{k}{L},\epsilon\right).
\]
The quantity $\frac{1}{L}\sum\limits_{k=1}^L\log\hat{c}\left(\frac{k}{L},\epsilon\right)$ represents a Riemann sum over $[0,1]$ for the fonction $\omega\mapsto\log\hat{c}(\omega,\epsilon)$. Thus an {\bf escape rate per lattice unit} can be defined in the limit of large periods 
\[
\lim\limits_{L\to\infty}\frac{\tilde\gamma_L(\epsilon)}{L}=\log\left(\frac{a}{2}\right)+\int_0^1\log\left(1-\epsilon(1-\cos2\pi\omega)\right)d\omega
\]
This normalized rate also decreases monotonically when $\epsilon$ increases. 

It is interesting to note that all the arguments above extend to CML with any convolution coupling. The expression for the escape rate per lattice unit remains unchanged
\[
\lim\limits_{L\to\infty}\frac{\tilde\gamma_L}{L}=\log\left(\frac{a}{2}\right)+\int_0^1\log|\hat{c}\left(\omega\right)|d\omega\quad\text{if}\quad a(2c_0-1)>2
\]
where $\hat{c}\left(\omega\right)=\sum\limits_{n\in\Z}c_ne^{2i\pi n\omega}$ is the (discrete) Fourier transform of the sequence $\{c_n\}$.

The explicit computations here crucially depend on the assumption that the derivative $(C\circ F_0)'$ is constant in phase space. In order to prove the existence of the escape rate per lattice unit for a broader class of CML, we use a more direct approach that relies on estimating the volume of the sets of points that remain located in the neighborhood of the repeller. 

\subsection{Approach via volumes}
A direct definition of the escape rate can be formulated as follows \cite{BR75,ER85}. Given a mapping $G:X\subset\R^n\to\R^n$ ($X$ bounded) with an invariant set $J$, the {\bf escape rate} from an $\eta$-neighborhood $U$ of $J$ ($\eta>0$) is the limit (if it exists) 
\[
\gamma(U)=-\lim\limits_{T\to\infty}\frac{1}{T}\log\text{Vol}(U_T)
\]
where $\text{Vol}(U_T)$ is the Lebesgue volume of the set $U_T$ of points $x$ such that $G^tx\in U$ for all $0\leq t<T$. In general, the limit does not depend on $U$ provided that $\eta$ is small enough. Moreover the escape rate of Axiom A attractors of diffeomorphisms on compact manifolds  satisfies the following variational principle \cite{BR75} 
\[
\gamma(U)=\sup\left\{\lambda^+_\mu-h_\mu\ :\ \mu\ \text{ergodic invariant probability, supp}(\mu)\subset J\right\}
\]
Eckmann and Ruelle \cite{ER85} conjectured that this relation should be valid in a more general setting. This has been recently proved for repellers of some one-dimensional maps \cite{BDM07} and a counter-example has also been constructed \cite{BBS99}. 

Back to CML on the lattice $\Z$, an escape rate can be introduced directly by limiting the dynamics to periodic configurations in $\Z_L$. Given an $\eta$-neighborhood $U\subset\ell^\infty(\Z)$ of the invariant set $K$ for the infinite CML, we consider the following limit 
\[
\gamma_L(U)=-\lim\limits_{T\to\infty}\frac{1}{T}\log\text{Vol}(U_T^L)
\] 
where $U_T^L$ is the set of periodic configurations $x\in\R^{\Z_L}$ such that $F^tx\in U$ for all $0\leq t\leq T$. (Of note, all configurations $F^tx$ must then also be $L$-periodic and the Lebesgue volume is computed in $\R^{\Z_L}$.)

\noindent
Following \cite{BS88,PS91}, we say that a convolution operator $C$ is of {\bf finite range} iff there exists $\zeta>1$ such that we have 
\[
\sum\limits_{n\in\Z}\zeta^{|n|}c_n<+\infty
\]
for the corresponding sequence $\{c_n\}_{n\in\Z}$. As announced before, the existence of the limit $\frac{\gamma_L(U)}{L}$ can be proved for any CML with short range couplings and local map satisfying the conditions (i)--(iii). 
\begin{Thm}
Given $f$ with properties {\rm (i)--(iii)} and a short range $C$ such that $\|\text{Id}-C\|<\epsilon_{f}$, for any sufficiently small $\eta$-neighborhood $U$ of $K$, the escape rate $\gamma_L(U)$ is well-defined for every $L>1$ and the following limit exists and is independent of $U$
\[
\gamma_\infty=\lim_{L\to\infty}\frac{\gamma_L(U)}{L}
\]
\label{THMTAUXFUIT}
\end{Thm}
The quantity $\gamma_\infty$ provides a definition of the (normalized) escape rate at the 'thermodynamics' limit in terms of estimates for the finite dimensional approximations. Of note, this definition does not rely on the existence of a conditionally invariant measure for the infinite system, nor on the direct computation of volumes in this limit. We conclude by a list of open questions
\begin{itemize}
\item[$\bullet$] Prove or disprove that  
\[
\gamma_\infty=\lim_{T\to\infty}\frac{1}{T}\lim_{L\to\infty}\frac{1}{L}\log\text{Vol}(U_T^L)
\]
(see comment after Lemma \ref{MAJOSUM} below). 
\item[$\bullet$] For the piecewise affine example of section 3.2, the quantities $\tilde\gamma_L$ and $\gamma_L(U)$ coincide (see Lemma \ref{BOUNDS} below). Prove that this relation holds for CML with arbitrary Pianigiani-Yorke maps.
\item[$\bullet$] Prove the Eckmann-Ruelle conjecture for periodic approximations of the CML.
\item[$\bullet$] Following \cite{B02}, prove that $\gamma_\infty$ corresponds to the escape rate of the $\Z^2$-action generated by $F$ and the spatial translations. 
\end{itemize}

\section{Proof of Theorem \ref{THMTAUXFUIT}}
\paragraph{First step: Reducing neighborhood to union of basic sets.}
We claim that w.l.o.g.\ one can rely on the set $\text{\bf I}\cap\R^{\Z_L}=\bigcup\limits_{\theta\in {\mathcal A}^{\Z_L}}(\text{\bf I}_\theta\cap\R^{\Z_L})$ \footnote{The equality here follows from the fact that the original intervals $\text{I}_i$ are disjoints; hence to be a periodic configuration in a basic set $\text{\bf I}_\theta$ imposes that the symbol configuration $\theta\in {\mathcal A}^{\Z_L}$ must be periodic.} to compute the escape rate of the periodic CML. 

\noindent
First, the set $K$ is included in the disjoint union $\text{\bf I}={\displaystyle\bigcup_{\theta\in {\mathcal 
A}^\Z}}\text{\bf I}_\theta$ and each component $K\cap \text{\bf I}_\theta$ is at positive distance 
from the complement $\ell^\infty(\Z)\setminus \text{\bf I}_\theta$ (see proof of condition (ii) in 
Appendix \ref{A-ERRATUM}). Therefore, there exists $\eta_0>0$ such that for every $\eta<\eta_0$, the 
$\eta$-neighborhood $U$ of $K$ 
is contained in $\text{\bf I}$. By passing to intersections 
with $\R^{\Z_L}$, it follows that the escape rate computed based on the $\eta$-neighborhood $U$ cannot exceed the rate computed with the union of basic sets. 

\noindent
To prove the converse inequality, consider an $\eta$-neighborhood $U$ contained in $\text{\bf I}$. The assumptions that the sets $\text{\bf I}_\theta$ are pairwise disjoints and that the intersection $F\text{\bf I}_\theta\cap \text{\bf I}_{\theta'}$ is non empty iff $\text{Int}\; F\text{\bf I}_\theta\supset \text{\bf I}_{\theta'}$ imply that
\[
\bigcap_{t=0}^{T} F^{-t}\text{\bf I}=\bigcup\limits_{\{\theta^t\}_{t=0}^{T}\in\Omega_{T+1}^{\Z}}\bigcap_{t=0}^{T}F^{-t}\text{\bf I}_{\theta^t}
\]
(where $\Omega_T$ represents the set of admissible words of length $T$ in ${\mathcal A}^{T}$). For every $T\in\N$, the latter set is contained in an $\eta_T$-neighborhood of $K$ for some $\eta_T>0$ such that$\lim\limits_{T\to\infty}\eta_T=0$. Therefore, for $T^*$ large enough such that $\eta_T<\eta$, we obtain that 
\[
\bigcap_{t=0}^{T^*-1} F^{-t}\text{\bf I}\subset U
\]
which implies that the escape rate obtained with the union of basic sets cannot exceed the rate computed based on neighborhoods. 

\noindent
Given a periodic symbol configuration $\theta\in {\mathcal A}^{\Z_L}$, let $\text{\bf I}^L_\theta=\text{\bf I}_\theta\cap\R^{\Z_L}$ be the corresponding set of periodic configurations. Analogous arguments to as before and the CML translation symmetry show that in order to estimate the escape rate of CML on periodic lattices, we may consider the volume of
\[
U_T^L=\bigcup\limits_{\{\theta^t\}_{t=0}^{T}\in\Omega_{T+1}^{\Z_L}}\bigcap\limits_{t=0}^{T}F^{-t}\text{\bf I}_{\theta^t}^L.
\]
\bigskip

\paragraph{Second step: Volume computation and partition function.}
The sets $\text{\bf I}_{\{\theta^t\}_{t=0}^{T}}^L=\bigcap\limits_{t=0}^{T}F^{-t}\text{\bf I}_{\theta^t}^L$ are pairwise disjoints; hence  
\[
\text{Vol}\left(U_T^L\right)=\sum\limits_{\{\theta^t\}_{t=0}^{T}\in\Omega_{T+1}^{\Z_L}}\text{Vol}(\text{\bf I}_{\{\theta^t\}_{t=0}^{T}}^L).
\]
The volumes in the RHS trivially read
\[
\text{Vol}(\text{\bf I}_{\{\theta^t\}_{t=0}^{T}}^L)=\int_{\text{\bf I}_{\{\theta^t\}_{t=0}^{T}}^L}dx
\]
The map $F^{T}$ is one-to-one on each $\text{\bf I}_{\{\theta^t\}_{t=0}^{T}}^L$ and we have $F\text{\bf I}^L_{\theta^t}\supset \text{\bf I}^L_{\theta^{t+1}}$ for all $0\leq t\leq T-1$ when $\{\theta^t\}_{t=0}^{T}$ is admissible. It results that $F^{T}\text{\bf I}_{\{\theta^t\}_{t=0}^{T}}^L= \text{\bf I}^L_{\theta^{T}}$. Introducing $y=F^Tx\in\R^{\Z_L}$, we obtain 
\begin{align*}
\text{Vol}\left(\text{\bf I}_{\{\theta^t\}_{t=0}^{T}}^L\right)=&\int_{\text{\bf I}_{\theta^{T}}^L}\prod_{t=0}^{T-1}\left|F'\circ F^{t-T}|_{\text{\bf I}_{\{\theta^k\}_{k=t}^T}}y\right|^{-1}dy\\
=&\left|C\right|_L^{-T}\int_{\text{\bf I}_{\theta^{T}}^L}\prod_{0\leq t<T\atop s\in\Z_L}\left|f'(F^{t-T}|_{\text{\bf I}_{\{\theta^k\}_{k=t}^T}}y)_s\right|^{-1}dy
\end{align*}
where $\left|C\right|_L$ is the determinant of the restriction $C\big{|}_{\R^{\Z_L}}$.

The last expression allows one to substitute the analysis of the asymptotic behaviour of $\text{Vol}\left(U_T^L\right)$ by the study of a kind of {\bf partition function} associated with the CML and defined by 
\[
Z_{L,T}=\sum\limits_{\{\theta^t\}_{t=0}^{T-1}\in\Omega_T^{\Z_L}}\sup\limits_{x\in \text{\bf I}_{\{\theta^t\}_{t=0}^{T-1}}^L}\prod_{0\leq t<T\atop s\in\Z_L}\left|f'(x_s^t)\right|^{-1}
\]
where we use the notation $x_s^t=(F^tx)_s$. The substitution is justified by the following result.
\begin{Lem}
There exist two numbers $0\leq\underline{c}<\overline{c}$ such that for all $L,T\geq 1$, we have  
\[
\underline{c}^L\left|C\right|_L^{-T}Z_{L,T}\leq \text{\rm Vol}\left(U_T^L\right) \leq \overline{c}^L\left|C\right|_L^{-T}Z_{L,T}
\]
\label{BOUNDS}
\end{Lem}
This statement leads to consider the quantity $K_L=\lim\limits_{T\to\infty}\frac{1}{T}\log Z_{L,T}$. If this limit exists, we have $\gamma_L(U)=K_L-\log\left|C\right|_L$ and, since 
\[
\lim\limits_{L\to\infty}\frac{\log\left|C\right|_L}{L}=\int_0^1\log|\hat{c}\left(\omega\right)|d\omega
\]
(see section 3.2), to obtain the theorem it only remains to show that the limit $\lim\limits_{L\to\infty}\frac{K_L}{L}$ exists.

\begin{proof}
The proof essentially consists in obtaining a bounded distorsion estimate (see for instance \cite{C96}), that is to say in proving the existence of a number $c_1\geq 1$ such that for any $T\geq1$ and every admissible word $\{\theta^t\}_{t=0}^{T-1}$, we have
\begin{equation}
c_1^{-1}\leq\frac{\prod_{t=0}^{T-1}\left|f'(x_s^t)\right|^{-1}}{\prod_{t=0}^{T-1}\left|f'(y_s^t)\right|^{-1}}\leq c_1,\quad\forall s\in\Z_L,\ x,y\in \text{\bf I}_{\{\theta^t\}_{t=0}^{T-1}}^L
\label{DISTORSION}
\end{equation}
To begin, notice that the definition of $F$ implies that when $C$ is invertible, the following inequality holds for every $L$-periodic configurations $x,y$
\[
\left|f(x_s)-f(y_s)\right|\leq \left\|C^{-1}\right\| \left\|Fx-Fy\right\|\quad \forall s\in\Z_L. 
\]
If, in addition, both configurations $x$ and $y$ belong to the same basic set $\text{\bf I}_\theta^L$, then we have 
\[
\left|x_s-y_s\right|\leq \alpha\left\|Fx-Fy\right\|\quad \forall s\in\Z_L
\]
where 
\[
\alpha=\frac{\left\|C^{-1}\right\|}{\inf\limits_{\text{I}}\left|f'\right|}\leq \frac{1}{(1-\epsilon_f)\inf\limits_{\text{I}}\left|f'\right|}<1
\]
due to the assumption $\|\text{Id}-C\|<\epsilon_{f}$ (see Appendix \ref{A-ERRATUM}). If the configurations $x$ and $y$ belong to the same cylinder set $\text{\bf I}_{\{\theta^t\}_{t=0}^{T-1}}^L$, then the argument can be repeated to obtain 
\[
\left|x^t_s-y^t_s\right|\leq \alpha^{T-t}\left\|x^T-y^T\right\|\quad \forall 0\leq t<T,\ s\in\Z_L
\]
In particular, we have $\left|x^t_s-y^t_s\right|\leq \alpha^{T-t}M$ where $M=2\max\limits_{\text{I}}|f|$.
Furthermore, $f$ is supposed to be of class $C^2$ on each interval $\text{I}_i$. It follows that  
\begin{equation}
\left|\log\left|f'(x)\right|-\log\left|f'(y)\right|\right|\leq\max\limits_{z\in \text{I}}\left|\frac{f''(z)}{f'(z)}\right| |x-y|\quad \forall x,y\in \text{I}_i
\label{ACCROIFIN}
\end{equation}
By combining the last two inequalities, we get 
\[
\left|\log\left|\frac{f'(x_s^t)}{f'(y_s^t)}\right|\right|\leq \beta M \alpha^{T-t}\quad \forall 0\leq t<T,\ s\in\Z_L,\ x,y\in \text{\bf I}_{\{\theta^t\}_{t=0}^{T-1}}^L
\]
where $\beta=\max\limits_{z\in \text{I}}\left|\frac{f''(z)}{f'(z)}\right|$. By summing over $t$ from 0 to $T-1$, the inequality (\ref{DISTORSION}) then easily follows with $c_1=\exp\frac{\alpha\beta M}{1-\alpha}$. 
\bigskip

Now, the inequality (\ref{DISTORSION}) implies that for all $x\in \text{\bf I}_{\{\theta^t\}_{t=0}^{T}}^L$, we have
\[
c_1^{-L}\sup_{x\in \text{\bf I}_{\{\theta^t\}_{t=0}^{T-1}}^L}\prod_{0\leq t<T\atop s\in\Z_L}\left|f'(x_s^t)\right|^{-1}\leq \prod_{0\leq t<T\atop s\in\Z_L}\left|f'(x_s^t)\right|^{-1},
\]
which yields  
\[
c_1^{-L}\left|C\right|_L^{-T}\text{Vol}(\text{\bf I}_{\theta^{T}}^L)\sup_{x\in \text{\bf I}_{\{\theta^t\}_{t=0}^{T-1}}^L}\prod_{0\leq t<T\atop s\in\Z_L}\left|f'(x_s^t)\right|^{-1}\leq \text{Vol}(\text{\bf I}_{\{\theta^t\}_{t=0}^{T}}^L).
\]
On the other hand, we obviously have 
\[
\prod_{0\leq t<T\atop s\in\Z_L}\left|f'(x_s^t)\right|^{-1}\leq \sup_{x\in \text{\bf I}_{\{\theta^t\}_{t=0}^{T-1}}^L}\prod_{0\leq t<T\atop s\in\Z_L}\left|f'(x_s^t)\right|^{-1},
\]
which implies
\[
\text{Vol}(\text{\bf I}_{\{\theta^t\}_{t=0}^{T}}^L)\leq \left|C\right|_L^{-T}\text{Vol}(\text{\bf I}_{\theta^{T}}^L)\sup_{x\in \text{\bf I}_{\{\theta^t\}_{t=0}^{T-1}}^L}\prod_{0\leq t<T\atop s\in\Z_L}\left|f'(x_s^t)\right|^{-1}
\]
The Lemma then immediately follows with $\underline{c}=\frac{1}{c_1}\min\limits_i\left|\text{I}_i\right|$ and $\overline{c}=\max\limits_i\left|\text{I}_i\right|$.
\end{proof}
\bigskip

\paragraph{Third step: Existence of the escape rate $\gamma_L(U)$.}
To begin, notice that the condition $x\in\text{\bf I}_{\{\theta^t\}_{t=0}^{T_1+T_2-1}}^L$ implies $x^{T_1}\in \text{\bf I}_{\{\theta^{t+T_1}\}_{t=0}^{T_2-1}}^L$. Therefore, for every triple $L,T_1,T_2\geq 1$, we have 
\begin{align*}
Z_{L,T_1+T_2}\leq&\sum\limits_{\{\theta^t\}_{t=0}^{T_1+T_2-1}\in\Omega_{T_1+T_2}^{\Z_L}}\left(\sup\limits_{x\in \text{\bf I}_{\{\theta^t\}_{t=0}^{T_1-1}}^L}\prod_{0\leq t<T_1\atop s\in\Z_L}\left|f'(x_s^t)\right|^{-1}\right.\\
&\left.\times\sup\limits_{x\in \text{\bf I}_{\{\theta^{t+T_1}\}_{t=0}^{T_2-1}}^L}\prod_{0\leq t<T_2\atop s\in\Z_L}\left|f'(x_s^{t+T_1})\right|^{-1}\right)\\
\leq&\sum\limits_{\{\theta^t\}_{t=0}^{T_1-1}\in\Omega_{T_1}^{\Z_L}}
\sup\limits_{x\in \text{\bf I}_{\{\theta^t\}_{t=0}^{T_1-1}}^L}
\prod_{0\leq t<T_1\atop s\in\Z_L}\left|f'(x_s^t)\right|^{-1}\\
&\times\sum\limits_{\{\theta^t\}_{t=0}^{T_2-1}\in\Omega_{T_2}^{\Z_L}}
\sup\limits_{x\in \text{\bf I}_{\{\theta^t\}_{t=0}^{T_2-1}}^L}
\prod_{0\leq t<T_2\atop s\in\Z_L}\left|f'(x_s^t)\right|^{-1}\\
=&Z_{L,T_1}Z_{L,T_2}
\end{align*} 
This property implies that the sequence $\{\log Z_{L,T}\}_{T\geq 1}$ is sub-additive; thus the following exists 
\[
\lim\limits_{T\to\infty}\frac{1}{T}\log Z_{L,T}=\inf\limits_{T\geq 1}\frac{1}{T}\log Z_{L,T}
\]
which, together with Lemma \ref{BOUNDS}, proves the existence of $\lim\limits_{T\to\infty}\frac{1}{T}\log\text{Vol}(U_L^T)$. 
\bigskip

\paragraph{Fourth step: Analysis of the asymptotic behaviour of $\gamma_L(U)$.}
Since for $\|\text{Id}-C\|<\epsilon_f$, the symbolic dynamics is given by the direct product $\Omega^\Z$, the spatial product in the definition of $Z_{L,T}$ can be decomposed into two independent spatial products, namely 
\begin{align*}
Z_{L_1+L_2,T}\leq&\sum\limits_{\{\theta^t_1\}_{t=0}^{T-1}\in\Omega_{T}^{\Z_{L_1}}}\sup\limits_{x\in \text{\bf I}_{\{(\theta_1\theta_2)^t\}_{t=0}^{T-1}}^{L_1+L_2}}\prod_{0\leq t<T\atop 1\leq s\leq L_1}\left|f'(x_s^t)\right|^{-1}\\
&\times\sum\limits_{\{\theta^t_2\}_{t=0}^{T-1}\in\Omega_{T}^{\Z_{L_2}}}\sup\limits_{x\in \text{\bf I}_{\{(\theta_1\theta_2)^t\}_{t=0}^{T-1}}^{L_1+L_2}}\prod_{0\leq t<T\atop L_1+1\leq s\leq L_1+L-2}\left|f'(x_s^t)\right|^{-1}
\end{align*}
where $(\theta_1\theta_2)^t$ denotes the spatial concatenation of the words $\theta_1^t$ et $\theta_2^t$ (which have spatial length $L_1$ and $L_2$ respectively). Accordingly, a relationship has to be established between the supremum involved in the first term of the RHS and the supremum in the definition of $Z_{L_1,T}$, and similarly for the second term with $Z_{L_2,T}$. To that goal we are going to use the following inequalities (which are consequences of (\ref{ACCROIFIN}))
\[
\exp-\beta\sum_{0\leq t<T\atop 1\leq s\leq L}\left|x_s^t-y_s^t\right|\leq
\prod_{0\leq t<T\atop 1\leq s\leq L}\frac{\left|f'(x_s^t)\right|}{\left|f'(y_s^t)\right|}\leq
\exp\beta\sum_{0\leq t<T\atop 1\leq s\leq L}\left|x_s^t-y_s^t\right|
\] 
where $x_s^t,y_s^t\in \text{I}_{\theta_s^t}$. It remains to control the double sum involved in these inequalities. This is the scope of the following statement whose proof is given in Appendix \ref{A-MAJOSUM}.
\begin{Lem}
For every $L,T\geq 1$ there exists $\Sigma_{L,T}>0$ with the following property
\[
\sup\limits_{L\geq 1}\lim\limits_{T\to\infty}\frac{1}{T}\Sigma_{L,T}<\infty
\]
such that for every admissible word $\{\theta^t\}_{t=0}^{T-1}\in\Omega_{T}^{\Z_L}$ and every configuration $x,y\in {\displaystyle\bigcap_{t=0}^{T-1}}F^{-t}\text{\bf I}$ with $x_s^t,y_s^t\in \text{I}_{\theta_s^t}$ for all $0\leq t<T$ and $1\leq s\leq L$, we have 
\[
\sum_{0\leq t<T\atop 1\leq s\leq L}\left|x_s^t-y_s^t\right|\leq \Sigma_{L,T}
\]
\label{MAJOSUM}
\end{Lem}
When combined with the upper bound above for $Z_{L_1+L_2,T}$, this Lemma implies the following inequality
\[
Z_{L_1+L_2,T}\leq Z_{L_1,T}Z_{L_2,T}e^{\beta(\Sigma_{L_1,T}+\Sigma_{L_2,T})}
\]
Consequently, the limit $\lim\limits_{T\to\infty}\frac{1}{T}\log Z_{L,T}$ is a sub-additive sequence of $L$ up to a constant. Proposition 9.6.4 in \cite{KH95} implies the existence of the double limit $\lim\limits\limits_{L\to\infty}\frac{1}{L}\lim\limits_{T\to\infty}\frac{1}{T}\log Z_{L,T}$, which concludes the proof of Theorem \ref{THMTAUXFUIT}. 

\noindent
Of note, the method of proof in Appendix \ref{A-MAJOSUM} does not allow to guarantee that  
\[
\lim\limits_{T\to\infty}\frac{1}{T}\lim\limits_{L\to\infty}\frac{1}{L}\Sigma_{L,T}=0
\]
from which commutation of space and time limits would follow in the definition of $\gamma_\infty$.

\appendix
\section{Proof of Pianigiani-Yorke conditions for weakly coupled map lattices}\label{A-ERRATUM}
We want to prove that, when $\|\text{Id}-C\|$ is sufficiently small (depending only on $f$), the symbolic dynamics remains unaffected and the conditions (i) and (ii) hold for the CML $C\circ F_0$ with infinite collection $\{\text{\bf I}_{\theta}\}_{\theta\in {\mathcal A}^\Z}$. The transitivity condition (iii) then naturally follows from the assumption that $f$ is transitive. 

The condition (ii) implies that the cylinder set $\bigcap\limits_{t=0}^{T} F^{-t}\text{\bf I}_{\theta^t}$ is non-empty for every admissible word $\{\theta^t\}_{t=0}^{T}\in\Omega_{T+1}^\Z$. The condition (i) guarantees that its diameter goes to 0 as $T\to\infty$. Hence $\bigcap\limits_{t=0}^{\infty} F^{-t}\text{\bf I}_{\theta^t}$ consists of a single point. The other properties in Theorem \ref{STABSTRUCT} are easily shown and left to the reader.
\bigskip

\paragraph{Proof of the condition (i).} The map $f$ being $C^2$ on $\text{I}$ implies that $F_0$ and thus $C\circ F_0$ are $C^2$ on $\text{\bf I}$. Moreover, if $\|\text{Id}-C\|\leq \epsilon<1$, then the inverse $C^{-1}$ exists and we have
\[
\left\|C^{-1}\right\|\leq \sum_{k\in\N}\|\text{Id}-C\|^k\leq\frac{1}{1-\epsilon},
\]
and thus $\|Cu\|\geq (1-\epsilon)\|u\|$ for all $u\in\ell^\infty(\Z)$. It follows that the norm $\|(C\circ F_0)'\|$ of the derivative is bounded below on $\text{\bf I}$ by $(1-\epsilon)\inf\limits_{\text{I}}\left|f'\right|$, a bound which is larger than 1 when $\epsilon$ is sufficiently small. 
\bigskip

\paragraph{Proof of the condition (ii).} As before, we assume that $\|\text{Id}-C\|\leq\epsilon<1$. 
We need to prove that, if $\epsilon$ is small enough, the cylinders $\text{\bf I}_\theta$ 
and their images satisfy the same dichotomy as in the uncoupled case: 
for any $\theta,\theta'\in {\mathcal A}^\Z$,

(a) either there exists $s\in \Z$ such that 
$\text{I}_{\theta'_s}\cap f\left(\text{I}_{\theta_s}\right)=\emptyset$, in which case
$$
\text{\bf I}_{\theta'}\cap F\text{\bf I}_\theta=\emptyset
$$

(b) or $\text{I}_{\theta'_s}\subset \text{Int}\; f\left(\text{I}_{\theta_s}\right)$ for all $s\in\Z$, in which case 
$$
\text{\bf I}_{\theta'}\subset  \text{Int}\; F\text{\bf I}_\theta
$$

Let us first denote, for any $i\in\mathcal{A}$, $f(I_{i})=[\alpha_{i}, \beta_{i}]$. 
We note also that, by compactness, we can quantify the properties of the local map $f$, {\it i.e.} 
choose $\gamma>0$ and $\delta>0$ such that, for any $i,j\in\mathcal{A}$,
\begin{align*}
I_{j}\cap f(I_{i})=\emptyset \quad 
& \Longrightarrow \quad I_{j}\cap [\alpha_{i}-\gamma,\beta_{i}+\gamma]=\emptyset\\
I_{j}\subset \text{Int}\; f\left(\text{I}_{\theta_s}\right) \quad 
& \Longrightarrow \quad I_{j}\subset (\alpha_{i}+\delta,\beta_{i}-\delta)
\end{align*}

Let also $M=2\max\limits_{\text{I}}|f|<\infty$. Then, for any $x\in F_0\text{\bf I}$ and any $s\in \Z$, we have
$$
|x_{s}-(C(x))_{s}|\leq \|\text{Id}-C\| \|x\|\leq \epsilon\frac{M}2.
$$
This estimate directly implies the following inclusion as soon as $\epsilon\leq \frac{2\gamma}M$
$$
F\text{\bf I}_{\theta}=C\left(\bigotimes_{s\in\Z}[\alpha_{\theta_s},\beta_{\theta_s}]\right)\subset
\bigotimes_{s\in\Z}[\alpha_{\theta_s}-\gamma,\beta_{\theta_s}+\gamma].
$$
From this, one immediately gets that $\text{\bf I}_{\theta'}\cap F\text{\bf I}_\theta=\emptyset$ 
in case (a).

For the other case, let us first remark that, by the method of point (i), one has 
$\|\text{Id}-C^{-1}\|\leq\frac\epsilon{1-\epsilon}$. 
Hence just as before, as soon as 
$\epsilon<\frac\delta{\delta+r} \Leftrightarrow\frac\epsilon{1-\varepsilon}<\frac\delta r$ 
and for any $\theta$ from case (b), one has
$$
C^{-1}\left(\bigotimes_{s\in\Z}[\alpha_{\theta_s}+\delta,\beta_{\theta_s}-\delta]\right)\subset
\bigotimes_{s\in\Z}[\alpha_{\theta_s},\beta_{\theta_s}],
$$
which implies 
$$
\text{\bf I}_{\theta'}\subset  \text{Int}\; \bigotimes_{s\in\Z}[\alpha_{\theta_s}
+\delta,\beta_{\theta_s}-\delta] \subset \text{Int}\; C \left(\bigotimes_{s\in\Z}
[\alpha_{\theta_s},\beta_{\theta_s}]\right)=\text{Int}\; F\text{\bf I}_\theta
$$ 

\section{Proof of Lemma \ref{MAJOSUM}}\label{A-MAJOSUM}
Let $\iota={\big(\inf\limits_{\text{I}}|f'|\big)}^{-1}$, let again $M$ be the diameter of the smallest ball containing $F\text{\bf I}$. We assume that we assume that $\|\text{Id}-C\|< 1$. As a consequence, the inverse $C^{-1}$ exists and is a convolution operator \cite{AF00}. Let $\{c_n^{(-1)}\}_{n\in\Z}$ be the sequence representing $C^{-1}$. If $x,y\in\R^\Z$ are two configurations such that $x_s,y_s\in \text{I}_{\theta_s}$ for all $1\leq s\leq L$, we have
\begin{align*}
|x_s-y_s|\leq &\iota\sum\limits_{n\in\Z}\left|c_n^{(-1)}\right|\left|\left(Fx\right)_{s-n}-\left(Fy\right)_{s-n}\right|\\
\leq &\iota M\sum\limits_{n<s-L}\left|c_n^{(-1)}\right|+\iota\sum\limits_{n=s-L}^{s-1}\left|c_n^{(-1)}\right|\left|\left(Fx\right)_{s-n}-\left(F_x\right)_{s-n}\right|\\
&+\iota M\sum\limits_{n>s-1}\left|c_n^{(-1)}\right|
\end{align*}
Let $P_{[1,L]}$ be defined in $\R^\Z$ as the canonical projection
\[
\left(P_{[1,L]}x\right)_s=\left\{\begin{array}{cl}
x_s&\text{if}\ 1\leq s\leq L\\
0&\text{otherwise}
\end{array}\right.
\]
and consider the configurations $\overline{H}$ et $\underline{H}$ defined by 
\[
\underline{H}_s=\left\{\begin{array}{cl}
0&\text{if}\ s\geq 0\\
M&\text{if}\ s< 0
\end{array}\right.
\quad\text{and}\quad
\overline{H}_s=\left\{\begin{array}{cl}
M&\text{if}\ s\geq 0\\
0&\text{if}\ s< 0
\end{array}\right.
\]
Consider also the operator ${\mathcal C}$ defined by 
\[
({\mathcal C}x)_s=\iota\sum\limits_{n\in\Z}\left|c_n^{(-1)}\right|x_{s-n},\quad \forall s\in\Z,\ x\in\ell^\infty(\Z)
\]
and the translation operator defined by $(Rx)_s=x_{s-1}$. The inequality above can be rewritten under the following compact form
\[
\left|x_s-y_s\right|\leq \left({\mathcal C}P_{[1,L]}\left(\left|Fx-Fy\right|\right)\right)_s+\left({\mathcal C}R\left(\underline{H}+R^{L}\overline{H}\right)\right)_s
\]

\noindent
By repeating the argument for the iterates, it results that if $x,y\in\R^\Z$ are two configurations such that $x_s^t,y_s^t\in \text{I}_{\theta_s^t}$ for all $1\leq s\leq L$ and $0\leq t<T$, then we have
\[
|x_s^t-y_s^t|\leq \left({\mathcal C}P_{[1,L]})^{T-t}\left(\underline{H}+\overline{H}\right)\right)_s+\left(\sum\limits_{k=0}^{T-t-1}\left({\mathcal C}P_{[1,L]}\right)^k{\mathcal C}R\left(\underline{H}+R^{L}\overline{H}\right)\right)_s
\]
It follows that $\sum\limits_{0\leq t<T \atop 1\leq s\leq L}|x_s^t-y_s^t|\leq \Sigma'_{L,T}+\Sigma''_{L,T}$ where
\[
\Sigma'_{L,T}=\sum\limits_{0\leq t<T \atop 1\leq s\leq L}\left(\left({\mathcal C}P_{[1,L]}\right)^{T-t}\left(\underline{H}+\overline{H}\right)\right)_s
\]
and
\[
\Sigma''_{L,T}=\sum\limits_{1\leq t\leq T \atop 1\leq s\leq L}t\left(\left({\mathcal C}P_{[1,L]}\right)^{T-t}{\mathcal C}R\left(\underline{H}+R^{L}\overline{H}\right)\right)_s
\]
In particular, this implies the existence of an upper bound for the finite sums $\sum\limits_{0\leq t<T \atop 1\leq s\leq L}|x_s^t-y_s^t|$. It remains to prove the asymptotic behaviours.

The condition $\|\text{Id}-C\|<\epsilon_{f}$ implies that $\|{\mathcal C}\|<1$ and thus $\left\|{\mathcal C}P_{[1,L]}\right\|<1$. Writing the Neumann series for this operator yields the following bounds (where the second inequality follows from the fact that the entries of $\left({\mathcal C}P_{[1,L]}\right)^t$ are all positive).
\[
\Sigma'_{L,T}\leq \sum_{t\geq 1 \atop 1\leq s\leq L}\left(\left({\mathcal C}P_{[1,L]}\right)^{t}\left(\underline{H}+\overline{H}\right)\right)_s\leq M L\sum_{t\geq 1}\|{\mathcal C}\|^t\leq\frac{M L\|{\mathcal C}\|}{1-\|{\mathcal C}\|}
\]

Let us now study the limit $\lim\limits_{T\to\infty}\frac{\Sigma''_{L,T}}{T}$. A direct calculation shows that, for every $x\in [0,1)$, the following limit holds 
\[
\lim\limits_{T\to\infty}\frac{1}{T}\sum\limits_{t=1}^{T}t x^{T-t}=\frac{1}{1-x}.
\]
This relation extends to operator components when the corresponding norm is smaller than 1. Consequently, using that $\left\|{\mathcal C}P_{[1,L]}\right\|<1$, we have
\begin{align*}
\lim\limits_{T\to\infty}\frac{1}{T}\sum\limits_{t=1}^{T}t\left(({\mathcal C}P_{[1,L]})^{T-t}{\mathcal C}R\left(\underline{H}+R^{L}\overline{H}\right)\right)_s=&\left|\left(\left(\text{Id}-{\mathcal C}P_{[1,L]})^{-1}{\mathcal C}R(\underline{H}\right.\right.\right.\\
&+\left.\left.\left. R^{L}\overline{H}\right)\right)_s\right|
\end{align*}
The assumption that $C$ is of finite range in Theorem \ref{THMTAUXFUIT} implies that the response $(C^{-1}\Delta)_n=c_n^{(-1)}$ (where $\Delta_s=\delta_{s,0}$) is exponentially localized around 0, {\sl i.e.}\ there exist $m_1>0$ and $0<\zeta_1<1$ such that \cite{BM97}
\[
\iota|c_n^{(-1)}|\leq m_1\zeta_1^{|n|},\quad \forall n\in\Z
\]
Using that any convolution operator commutes with pointwise limits of boun\-ded sequences in $\ell^\infty(\Z)$ \cite{AF00}, it follows that  
\[
({\mathcal C}R\underline{H})_s\leq Mm_1\sum\limits_{n\geq s}\zeta_1^{|n|}\quad\text{and}\quad
({\mathcal C}R^{L+1}\overline{H})_s\leq Mm_1\sum\limits_{n\leq s-L-1}\zeta_1^{|n|},\quad \forall s\in\Z
\]

\noindent
Furthermore, the operator $\text{Id}-{\mathcal C}P_{[1,L]}$ is accordingly of finite range for every $L>1$ (when the finite range property is extended to bounded linear operators $A$ of the form
\[
(A x)_s=\sum\limits_{n\in\Z}a_{n,s}x_{s-n},\quad \forall s\in\Z
\]
{\sl i.e.}\ $\sup\limits_{s\in\Z}\sum\limits_{n\in\Z}\zeta^{|n|}|a_{n,s}|<\infty$ for some $\zeta>1$.) The associated decay rate $\zeta_{C}$ is identical to that of $C^{-1}$ and is thus independent of $L$. Using again \cite{BM97}, it results that, for every $n\in\Z$, the response $\left(\text{Id}-{\mathcal C}P_{[1,L]}\right)^{-1}R^n\Delta$ is exponentially localized around $n$. Moreover the decay rate is independent of $L$ because it is controlled by $\zeta_{C}$ and $\frac1{1-\|{\mathcal C}\|}$). Hence there exist $m_2>0$ and $0<\zeta_2<1$ such that for all $n\in\Z$, we have
\[
\left|\left(\left(\text{Id}-{\mathcal C}P_{[1,L]}\right)^{-1}R^n\Delta\right)_s\right|\leq m_2\zeta_2^{|n-s|},\quad\forall s\in\Z,L>1.
\]
Furthermore by using Neumann series and the decomposition $x=\sum\limits_{s\in\Z}x_sR^s\Delta$ for a configuration $x=\left\{x_s\right\}_{s\in\Z}\in\ell^\infty(\Z)$, one shows that 
\[
\left(\left(\text{Id}-{\mathcal C}P_{[1,L]}\right)^{-1}{\mathcal C}R\underline{H}\right)_s=\sum\limits_{n\in\Z}({\mathcal C}R\underline{H})_{n}\left(\left(\text{Id}-{\mathcal C}P_{[1,L]}\right)^{-1}R^n\Delta\right)_s,\quad \forall s\in\Z
\]

\noindent
Combining the previous estimates leads to 
\begin{align*}
\left|\left(\left(\text{Id}-{\mathcal C}P_{[1,L]}\right)^{-1}{\mathcal C}R\underline{H}\right)_s\right|\leq&Mm_1m_2\sum\limits_{n\in\Z}\zeta_2^{|n|}\sum\limits_{m\geq s-n}\zeta_1^{|m|}\\
\leq&\frac{Mm_1m_2}{1-\zeta_1}\left(\frac{\zeta_1^s}{1-\zeta_1\zeta_2}+\zeta_2^{s}(\frac{1+\zeta_1}{1-\zeta_2}+\frac{\zeta_1^2}{1-\zeta_1\zeta_2})\right.\\
&+\left.\zeta_2\frac{\zeta_1^s-\zeta_2^{s}}{\zeta_1-\zeta_2}\right),
\end{align*}
for all $s\geq 1$; and similarly
\begin{align*}
\left|\left(\left(\text{Id}-{\mathcal C}P_{[1,L]}\right)^{-1}{\mathcal C}R^{L+1}\overline{H}\right)_s\right|&\leq\frac{Mm_1m_2}{1-\zeta_1}\left(\frac{\zeta_1^{L+1-s}}{1-\zeta_1\zeta_2}\right.\\
&+\left.\zeta_2^{L+2-s}(\frac{1+\zeta_1}{1-\zeta_2}+\frac{\zeta_1^2}{1-\zeta_1\zeta_2})\right.\\
&+\left.\frac{\zeta_1^{L+2-s}-\zeta_2^{L+2-s}}{\zeta_1-\zeta_2}\right),
\end{align*}
for all $s\leq L$. Together with the upper bound above on $\Sigma'_{L,T}$, by summing over $s$, we finally conclude that 
\[
\sup\limits_{L\geq 1}\lim\limits_{T\to\infty}\frac{1}{T}\Sigma_{L,T}<\infty.
\]

\noindent
{\bf Acknowledgments}

\noindent
BF would like to thank the Courant Institute (New York University) for hospitality during the paper preparation. His support in this phase was partly provided by CNRS and by the EU Marie Curie fellowship PIOF-GA-2009-235741.

\end{document}